\documentclass[12pt]{amsart}

\usepackage{graphicx}

\usepackage[all,cmtip]{xy}

\usepackage{marvosym}

\usepackage{answers}
\Newassociation{solution}{Solution}{ans}

\theoremstyle{definition}
\newtheorem{definition}{Definition}
\newtheorem{example}{Example}
\newtheorem{remark}{Remark}

\theoremstyle{plain}

\newtheorem{theorem}{Theorem}

\author{J. R. Arteaga, M. Malakhaltsev}
\title[$G$-structures and differential equations]{Ideas of E.~Cartan and S.~Lie in modern geometry: $G$-structures and differential equations. 
Lecture 4}
\address{}
\email{}

\begin{document}
\maketitle
\Opensolutionfile{ans}[answers_lecture_4]
\medskip
\fbox{\fbox{\parbox{5.5in}{
\textbf{Problem:}\\
How to find invariants of singularities of a $G$-structure?
}}}
\vspace{1cm}

\section*{Simple example} 
Let $V = (v^1,v^2)$ be a vector field on $\mathbb{R}^2$.
If $V(x_0) \ne 0$, then $V$ restricted to a neighborhood $U(x_0)$ is equivalent to $\partial_1$, hence all nonvanishing vector fields are locally equivalent and there are no invariants.

If $V(x_0) = 0$, and $W(x_0) = 0$,  and the matrix $||\partial_i V^j(x_0)||$ is not similar to the matrix $||\partial_i W^j(x_0)||$, then $V$ and $W$ are not locally equivalent and there arise invariants.

Generally, \emph{non-regular points have their own invariants}.

\section{Jets of smooth maps}
Let $M$ and $N$ be smooth manifolds, $x \in M$ is a point. 
We will denote by $f : (M,x) \to N$ a smooth map which is defined in an open neighborhood $U$ of $x$. 
By $D(f)$ we will denote the domain $U$. 

Denote by $C^\infty((M,x),N)$ the set of all smooth maps $f : (M,x) \to N$.
On the set $C^\infty((M,x),N)$ we introduce the equivalence relation: we say that $f$, $g$ in $C^\infty( (M,x),N)$ are equivalent at $x \in M$ ($f \sim_x g$) if 
\begin{equation}
\exists W \in \mathcal{U}_x, \quad W \subset D(f) \cap D(g) \text{ such that } f|_W = g|_W.
 \label{eq:equivalence_between_local_maps}
\end{equation}

\begin{definition}
\label{def:germ_of_map}
A \emph{germ of a map} $f \in C^\infty( (M,x) ,N)$  at the point $x$ is the equivalence class of
$f \in C^\infty((M,x),N)$ with respect to $\sim_x$ denoted by $\langle f \rangle_x$. 
\end{definition}
We set 
\begin{multline}
\mathcal{G}_x(M,N) = C^\infty( (M,x) ,N) / \sim_x = 
\\
=
\left\{ \langle f \rangle_x \mid f \in C^\infty( (M,x) ,N) \right\}. 
\label{eq:set_of_germs}
\end{multline}

If $M_1$, $M_2$, and $M_3$ are smooth manifolds, and $f_1 \in C^\infty((M_1, x_1),M_2)$, $f_2 \in C^\infty((M_2,x_2),M_3)$ with $f(x_1) = x_2$, then we can define the composition of germs as follows:
\begin{eqnarray}
\langle f_2 \rangle_{x_2} \circ \langle f_1 \rangle_{x_1} = \langle f_2 \circ f_1 \rangle_{x_1}
\label{eq:comp_of_germs}
\end{eqnarray} 

Another equivalence relation $\sim_k$ on the set $C^\infty( (M,x),N)$ is introduced as follows:
for $f,g \in C^\infty_x(M,N)$ such that $f(x)=g(x)$ we take the coordinate systems $(U,x^i)$ in a neighborhood $U$ of $x$ and $(V,y^\alpha)$ in a neighborhood $V$ of $y$. 
Then we say that $f$ and $g$ are equivalent ($f \sim_k g$) if the Taylor series of the coordinate representations of $f$ and $g$ coincide up to the order $k$. One can prove that this equivalence relation does not depend on the choice of coordinate systems.  
\begin{definition}
The equivalence class $j^k_x f$ of $f$ is called \emph{$k$-jet of the map $f$ at the point $x$}.    
\end{definition}
The set of all $k$-jets of maps $f \in C^\infty( (M,x) ,N)$ will be denoted by $J^k_x(M,N)$.

It is clear that if $ f, g \in C^\infty( (M,x) ,N)$ determine the same germ at $x$, that is 
if $f \sim_x g$, then $j^k_x f = j^k_x g$ for any $k$. 
Also, the composition of maps, or the composition of germs, define the composition of $k$-jets:
\begin{equation}
j^k_{f_1(x)} (f_2) \circ j^k_x(f_1) = j^k_x( f_2 \circ f_1).
\label{eq:composition_of_jets}
\end{equation}

For the set of all $k$-jets $J^k(M,N)$ one can define two natural projections:
\begin{eqnarray}
\pi_0 : J^k(M,N) \to M, \quad j^k_{x}f \mapsto x,
\\
\pi_1 : J^k(M,N) \to N, \quad j^k_{x}f \mapsto f(x).
\label{eq:projections_on_jet_space}
\end{eqnarray}
The set $J^k(M,N)$, $\dim M = m$, $\dim N = n$ can be endowed by a manifold structure so that these projections are bundle projections.
We set 
\begin{equation}
J^k_{x,y}(M,N) = \left\{ j^k f \in J^k(M,N) \mid \pi_0(j^k f) = x, \pi_1(j^k f) = y \right\}. 
\label{eq:set_of_jets_joining_two_points}
\end{equation}
The typical fiber of the bundle $(J^k(M,N),\pi_0,M)$ is the manifold 
$J^k_{0,0}(\mathbb{R}^m,\mathbb{R}^n) \times N$, 
and of the bundle $(J^k(M,N),\pi_1,N)$ is $J^k_{0,0}(\mathbb{R}^m,\mathbb{R}^n) \times M$.

The manifold $J^k(M,M)$ endowed with the projections $\pi_0$ and $\pi_1$, and the composition of $k$-jets is a groupoid. 

If $\pi : E \to M$ is a fiber bundle, then the set of $k$-jets of sections of the bundle $E$ is denoted by $J^k(E)$.

\section{Bundle of germs of diffeomorphisms. Coframe bundle of $k$th order}

\subsection{Bundle of germs of diffeomorphisms}
Let $\mathcal{D}(m)$ be the group of all germs of local diffeomorphisms of $\mathbb{R}^m$ at $0 \in \mathbb{R}^m$:
\begin{equation*}
\mathcal{D}(m) = \left\{ \langle \varphi\rangle_0  \mid 
\varphi : (\mathbb{R}^m,0) \to (\mathbb{R}^m,0) 
\text{ is a local diffeomorphism} \right\} 
\end{equation*}
endowed with the operation of composition of germs: 
$\langle \varphi_1 \rangle_0 \circ \langle \varphi_2 \rangle_0 = 
\langle \varphi_1 \circ \varphi_2 \rangle_0$. 

Now let us consider 
\begin{equation*}
\mathcal{B}(M) = \left\{ \langle f \rangle_x \mid f : (M,x) \to (\mathbb{R}^n,0) 
\text{ is a local diffeomorphism }\right\}
\end{equation*} 
We have natural projection 
\begin{equation*}
\pi : \mathcal{B}(M) \to M, \quad \pi(\langle f\rangle_x)=x, 
\end{equation*}
and the natural right action of $\mathcal{D}(m)$ on $\mathcal{B}(M)$:
\begin{equation*}
\langle f\rangle_x \cdot \langle \varphi\rangle_0 = \langle \varphi^{-1} \circ f\rangle_x. 
\end{equation*}
so $(\mathbb{B}(M),\pi,M)$ can be considered as a ``principal fiber bundle''.  

Let $(U, u : U \to V \subset \mathbb{R}^m)$ be a coordinate map. 
This map determines the ``trivialization'' of the bundle $\mathcal{B}(M)$ over $U$. 
Let 
\begin{equation*}
t_a : \mathbb{R}^m \to \mathbb{R}^m, \quad t_a(v) = v + a
\end{equation*}
be the parallel translation of $\mathbb{R}^m$. 
Then  
\begin{equation*}
\mathcal{U}: \pi^{-1}(U) \to U \times \mathcal{D}(m), 
\quad
\langle f \rangle_x \to \left(x, \langle f \circ u^{-1} \circ t_{u(x)} \rangle_0\right).
\end{equation*} 
gives us the required trivialization.
The inverse map is 
\begin{equation*}
\mathcal{U}^{-1} : U \times \mathcal{D}(m) \to \pi^{-1}(U)  
\quad
(x, \langle \varphi \rangle_0) \to \langle t_{-u(x)} \circ u \circ f^{-1} \rangle.
\end{equation*}
 
Now assume that we have two coordinate systems $(U,u)$ and $(U,\bar u)$ on $M$. 
Then, 
\begin{equation*}
\bar{\mathcal{U}} \circ \mathcal{U}^{-1} : U \times \mathcal{D}(m) \to U \times \mathcal{D}(m),
\quad (x, \varphi) \to (x, t_{-\bar u(x)} \circ \bar u \circ u^{-1} \circ t_{u(x)}).
\end{equation*}
so the gluing functions of the atlas of the ``principal bundle'' $(\mathcal{B}(M),\pi,M)$ constructed by an atlas $(U_\alpha,u_\alpha)$ of the manifold $M$ are
\begin{equation*}
g_{\beta\alpha} : U_{\alpha} \cap U_\beta \to \mathcal{D}(m), \quad
 g_{\beta\alpha}(x) =  t_{-u_\beta(x)} \circ u_\beta \circ u^{-1}_\alpha \circ t_{u_\alpha(x)}.
\end{equation*} 

\begin{remark}
In what follows we will use unordered multiindices.
We denote by $\mathcal{I}(m)$ the set of all unordered multiindices 
$I = \left\{ i_1 i_2 \dots i_k \right\}$, where $1 \le i_l \le m$, for all $l=\overline{1,k}$. 
The number $k$ is called the length of the multiindex and is denoted by $|I|$. 
Also, we set $I_k(m) = \left\{ I \in \mathcal{I}(m) \mid |I| = k  \right\}$.
\end{remark}

\subsection{Differential group of $k$th order}

The $k$th order \emph{differential group} is the set of $k$-jets:  
\begin{equation*}
D^k(m) = \left\{ j^k_0(\varphi) \mid 
\varphi : (\mathbb{R}^m,0) \to (\mathbb{R}^m,0) 
\text{ is a local diffeomorphism } \right\}. 
\end{equation*}

On the set $D^k(m)$ consider the operation 
\begin{equation}
D^k(m) \times D^k(m) \to D^k(m), \quad j^k_0(\varphi) \cdot j^k_0(\psi) = j^k(\varphi \circ \psi),
\label{eq:operation_in_kth_order_differential_group}
\end{equation}
then $(D^k(m),\cdot)$ is a group.

Denote $\varphi^k_I = \left.\partial_I\right|_0 \varphi^k$. 
Then 
\begin{equation}
\mathcal{C}^k : D^k(m) \to \mathbb{R}^N, \quad  j^k_0(\varphi) \to \{\varphi^k_I\} 
\label{eq:coordinates_on_kth_order_differential_group}
\end{equation} 
in a one-to-one map of $D^k(m)$ onto the open set in $\mathbb{R}^N$ determined by the inequality 
$\det\|\varphi^k_i\| \ne 0$. 
In this way we get globally defined coordinates on $D^k(m)$ which will be called \emph{natural coordinates}. 
With respect to the natural coordinates the product \eqref{eq:operation_in_kth_order_differential_group} is written in terms of polynomials, therefore is a smooth map. 
Thus $D^k(m)$ is a Lie group. 

\subsection{Bundle of $k$th order holonomic coframes}
For an $m$-dimensional manifold $M$ consider the set
\begin{equation*}
B^k(M) = \left\{ j^k_x{f} \mid f : (M,x) \to (\mathbb{R}^m,0) 
\text{ is a local diffeomorphism }\right\}
\end{equation*} 
whose elements are called \emph{$k$-coframes} or \emph{coframes of order $k$} of the manifold $M$. 
We have the projection
\begin{equation*}
\pi^k : B^k(M) \to M, \quad \pi(j_x(f))=x. 
\end{equation*}

On the set $B^k(M)$ we have the right $D^k(m)$-action: 
\begin{equation*}
j^k_x(f) \cdot j^k_0(\varphi) =  j^k_0(\varphi^{-1} \circ f).
\end{equation*}
and one can easily prove that this action is free and its orbits are the fibers of the projection $\pi$. 

\subsubsection{Trivializing charts of $B^k(M)$. Gluing maps}
In what follows we set $t_a : \mathbb{R}^m \to \mathbb{R}^m$, $t_a(x) = x + a$, the parallel translation of $\mathbb{R}^m$ with respect to $a \in \mathbb{R}^m$.

Let $(U, u : U \to V \subset \mathbb{R}^m)$ be a coordinate chart on $M$. 
We have the one-to-one map
\begin{multline}
\mathcal{T}^k: (\pi^k)^{-1}(U) \to U \times D^k(m), 
\\
j^k_x(f) \to \left(x, j^k_0(t_{-u(x)} \circ u \circ f^{-1}) \right).
\label{eq:trivialization_of_kth_order_coframe_bundle}
\end{multline}
The map $\mathcal{T}^k$ is $D^k(m)$-equivariant because
\begin{multline*}
\mathcal{T}^k( j^k_x f \cdot j^k_0\varphi) = \mathcal{T}^k \left(j^k_x (\varphi^{-1} \circ f)\right) = 
 \left(x, j^k_0(t_{-u(x)} \circ u \circ f^{-1} \circ \varphi)\right)= 
\\
=\left(x, j^k_0(t_{-u(x)} \circ u \circ f^{-1}) \cdot j^k_0\varphi \right)= 
\left(x, j^k_0(t_{-u(x)} \circ u \circ f^{-1}) \right) \cdot j^k_0\varphi. 
\end{multline*}
Since $D^k(m)$ is a Lie group, the map $\mathcal{T}^k$ defines a trivializing chart for the map $\pi^k : B^k(M) \to M$.

Therefore, for each atlas $(U_\alpha,u_\alpha)$, we construct the atlas of trivializing charts 
$(U_\alpha,\mathcal{T}^k_\alpha)$.  Find the gluing maps for this atlas. 

Assume that $(U_\alpha,u_\alpha)$, $(U_\beta,u_\beta)$ are two coordinate systems on $M$, and 
$U_\alpha \cap U_\beta \ne \emptyset$.
Then 
\begin{equation*}
j^k_0(t_{-u_\beta(x)} \circ u_\beta \circ f^{-1}) =
j^k_0(t_{-u_\beta(x)} \circ u_\beta \circ u_\alpha^{-1} \circ t_{u_\alpha(x)}) \cdot 
j^k_0(t_{-u_\alpha(x)} \circ u_\alpha \circ f^{-1}) 
\end{equation*} 
Therefore, the gluing maps are 
\begin{equation}
g_{\beta\alpha} : U_\alpha \cap U_\beta \to D^k(m), \quad
g_{\beta\alpha}(x) =  
j^k_0(t_{-u_\beta(x)} \circ u_\beta \circ u_\alpha^{-1} \circ t_{u_\alpha(x)}) 
\label{eq:gluing_maps_for_B_k}
\end{equation} 

Since the gluing functions are smooth, we conclude that $\pi^k : B^k(M) \to M$ is a $D^k(m)$-principal bundle over $M$ which is called \emph{the bundle of $k$-coframes of $M$} or \emph{the bundle of coframes of order $k$ of $M$}.

For a coordinate chart $(U,u)$ there is defined a section
\begin{equation}
s : U \to B^k(M), \quad s(x) = j^k_x u, 
\label{eq:natural_sections_of_B_k}
\end{equation}
which is called the \emph{natural $k$-coframe field} associated with a coordinate chart $(U_\alpha,u_\alpha)$. 

\subsubsection{Natural coordinates on $B^k(M)$}
If $(U,u)$ is a coordinate chart on $M$, and $\mathcal{T}^k$ is the corresponding trivialization of $B^k(M)$. 
Then 
\begin{equation}
(u \times \mathcal{C}^k) \circ \mathcal{T}^k : (\pi^k)^{-1}(U) \to \mathbb{R}^m \times \mathbb{R}^N
\end{equation}
gives \emph{natural local coordinates on} $B^k(M)$.

The section $s : U \to B^k(M)$ \eqref{eq:natural_sections_of_B_k} is written with respect to this coordinate system as follows:
\begin{equation}
s(u^k) = (u^k, \delta^k_i, 0).
\label{eq:natural_section_wrt_natural_coordinates}
\end{equation} 

\begin{remark}
We have the natural projections $\pi^k_l : B^k(M) \to B^l(M)$, $k \ge l$, which are, in turn, principal fiber bundles with the group $H^k_l$ which is the kernel of the natural homomorphism $D^k(m) \to D^l(m)$. 
\end{remark}

\subsection{Case $k=1$}
The Lie group $D^1(m) \cong GL(m)$ and $B^1(M)=B(M)$ is the coframe bundle of $M$.

\subsection{Case $k=2$}
\subsubsection{The group $D^2(m)$}
Elements of the group $D^2(m)$ are the $2$-jets of germs $\varphi \in \mathcal{D}(m)$. The coordinate system 
\eqref{eq:coordinates_on_kth_order_differential_group} in this case is 
\begin{equation}
j^2_0\varphi \longrightarrow (\varphi^k_i, \varphi^k_{ij}), \text{  where } 
\varphi^k_i = \frac{\partial \varphi^k}{\partial u^i}(0),
\quad \varphi^k_{ij} = \frac{\partial^2 \varphi^k}{\partial u^i \partial u^j}(0).  
\end{equation}
Here $u^i$ are coordinates on $\mathbb{R}^m$, and it is clear that $\varphi^k_{ij} = \varphi^k_{ji}$. 
From this follows that $\dim D^2(m) = m^2 + m^2(m+1)/2$.

Now, if  
$j^2_0\varphi \rightarrow (\varphi^k_i, \varphi^k_{ij})$, $j^2_0\psi \rightarrow (\psi^k_i, \psi^k_{ij})$, and
\begin{equation*}
j^2_0\psi \cdot j^2_0\varphi = j^2_0(\psi\circ\varphi) \rightarrow  (\eta^k_i, \eta^k_{ij}),
\end{equation*}
by the chain rule we get that 
\begin{equation}
\eta^k_i = \psi^k_s \varphi^s_i, \quad \eta^k_{ij} = \psi^k_{pq} \varphi^p_i \varphi^q_j + \psi^k_s \varphi^s_{ij}
\label{eq:product_D_2}
\end{equation}  
These formulas express the product in the group $D^2(m)$ in terms of the natural coordinates $(\varphi^k_i,\varphi^k_{ij})$. 

\subsubsection{The bundle $B^2(M)$}
The elements of $B^2(M)$ are the $2$-jets of local diffeomorphisms $f : (M,x) \to (\mathbb{R}^m,0)$. 
The natural coordinates \eqref{eq:natural_section_wrt_natural_coordinates} in this case can be found as follows.
Let $(U,u)$ be a coordinate chart on $M$. Then,  for any $j^2_x f$ with $x \in U$, the diffeomorphism 
\begin{equation*}
f \circ u^{-1} : (\mathbb{R}^m, u(x)) \to (\mathbb{R}^m,0)
\end{equation*}
can be written as 
$w^k=f^k(u^i)$, where $w^k$ are standard coordinates on $\mathbb{R}^m$, and $f^k(u^i(x)) = 0$. 
Then take the inverse diffeomorphism $u^k = \widetilde{f}^k(w^i)$, and the local diffeomorphism $t_{-u(x)} \circ u \circ f^{-1}$ has the form $\widetilde{f}^k(w^i)-u^k(x)$.
Therefore, the natural coordinates of $j^k_x f$ induced by a coordinate chart $(U,u)$ on $M$ are 
$(u^k,u^k_i,u^k_{ij})$, where 
\begin{equation}
u^k_i = \frac{\partial\widetilde{f}^k}{\partial w^i}(0), \quad 
u^k_{ij} = \frac{\partial^2\widetilde{f}^k}{\partial w^i \partial w^j}(0). 
\label{eq:natural_coordinates_B_2}
\end{equation}   
The derivatives of $\widetilde{f}$ at $0$ can be expressed in terms of the derivatives of $f$ at $u^k(x)$. If we denote 
\begin{equation*}
f^k_i = \frac{\partial f^k}{\partial u^i}(u(x)),
f^k_{ij} = \frac{\partial^2 f^k}{\partial u^i \partial u^j}(u(x)). 
\end{equation*} 
then
\begin{equation*}
u^k_i = \widetilde{f}^k_i, \quad u^k_{ij} = - \widetilde{f}^k_s f^s_{lm} \widetilde{f}^l_i \widetilde{f}^m_j.
\end{equation*} 

With respect to the natural coordinate system the $D^2(m)$-action is written as follows:
\begin{equation}
(u^k, u^k_i, u^k_{ij}) \cdot (\varphi^k_i, \varphi^k_{ij}) = (u^k, u^k_s \varphi^s_i,   
u^k_{pq} \varphi^p_i \varphi^q_j + u^k_s \varphi^s_{ij}).
\label{eq:D_2(m)-action_wrt_natural_coordinates}
\end{equation}

Let us express the gluing maps in terms of the natural coordinates. If $(U,u)$ and $(U',u')$ are coordinate charts on $M$ sucht that $U \cap U' \ne \emptyset$, then from \eqref{eq:gluing_maps_for_B_k} it follows that the corresponding gluing map is  
\begin{equation*}
g : U \cap U' \to D^2(m), g(x) = j^2_0(t_{-u'(x)} \circ u' \circ u^{-1} \circ t_{u(x)}) 
\end{equation*}  
If the coordinate change $u' \circ u^{-1}$ is written as $v^k = v^{k}(u^i)$, then we have to take derivatives at $0$ of the map $v^{k}(u^i+u^i(x))-v^k(u^i(x))$, which are equal to the derivatives of the functions $v^k$ at $u^k(x)$. Therefore,
\begin{equation}
g : U \cap U' \to D^2(m), \quad 
g(x) = 
\left(\frac{\partial v^k}{\partial u^i}(u(x)),  \frac{\partial v^k}{\partial u^i \partial u^j}(u(x))\right).
\label{eq:gluing_maps_for_B_2}
\end{equation} 

Therefore, \emph{the bundle $B^2(M) \to M$ is the $D^2(m)$-principal bundle with gluing maps \eqref{eq:gluing_maps_for_B_2}}.

\section{First prolongation of a $G$-structure}
\subsection{First prolongation of an integrable $G$-structure}
Let $P(M) \to M$ be an integrable $G$-structure, that is a subbundle of $B(M)$ such that there exists an atlas 
$\mathcal{A} = (U_\alpha,u_\alpha)$ such that the natural coframes of the atlas are sections of $P(M)$, or equivalently, 
the coordinate change $u^{k'} = u^{k'}(u^i)$ has the property that $\|\frac{\partial u^{k'}}{\partial u^k}\| \in G$.  

In this case, we can specify the set $\mathcal{B}_G$ of local diffeomorphisms  $f : (M,x) \to (\mathbb{R}^m,0)$ such that for each coordinate map $u$ of the atlas $\mathcal{A}$, the local diffeomorphism $f \circ u^{-1}$ has the Jacobi matrix at in $G$ at all points of its domain. 
It is clear that 
\begin{equation}
P(M) = \left\{ j^1_x f  \mid \| \frac{\partial (f \circ u^{-1})^k}{\partial u^i}|_{u(x)} \| \in G  \right\}
\end{equation}
and consider 
\begin{equation}
P^1(M) = \left\{ j^2_x f  \mid f \in \mathcal{B}_G \right\}
\label{eq:first_holonomic_prolongation_of_P}
\end{equation}
\begin{remark}
Note that if $j^2_x f$ is an element of $P^1(M)$, then the Jacobi matrix 
$\|\left.\frac{\partial (f\circ u^{-1})^k}{\partial u^i}\right|_{u(x)}\|$ 
is an element of $G$. 
However, the converse is not true because, by definition, the Jacobi matrix belongs to $G$ for each point which implies conditions on the second derivatives of $f \circ u^{-1}$ (see below). 
\end{remark}

In the same manner introduce the set $\mathcal{D}_G$ of local diffeomorphisms  
$\varphi : (\mathbb{R}^m,0) \to (\mathbb{R}^m,0)$ whose Jacobi matrices are elements of $G$ at all points of their domains. 
Consider the Lie subgroup of $D^2(m)$:
\begin{equation}
G^1 = \left\{ j^2_0 \varphi \in D^2(m) \mid \varphi \in \mathcal{D}_G \right\}
\label{eq:first_holonomic_prolongation_of_G}
\end{equation}
The Lie subgroup $G^1 \subset D^2(m)$ is called the \emph{first holonomic prolongation of the group} $G$. 

The subset $P^1(M)$ is a submanifold of $B^2(M)$, and is the total space of a principal subbundle of $B^2(M) \to M$ with the subgroup $G^1 \subset D^2(m)$. This subbundle is called  
the  \emph{first holonomic prolongation of the integrable $G$-structure $P$}. 

It is also clear that, for the atlas $\mathcal{A}$ the gluing map \eqref{eq:gluing_maps_for_B_2} takes values in the subgroup $G^1$. Therefore, \emph{an integrable $G$-structure defines reduction of the principal bundle $B^2(M)$ to the structure group $G^1 \subset D^2(m)$}. 

\subsection{Algebraic structure of the Lie group $G^1$} 
We have the surjective Lie group morphism 
\begin{equation}
p^1 : G^1 \to G, \quad j^2_0 \varphi \mapsto j^1_0 \varphi. 
\label{eq:surjective_morphism_G_1_to_G}
\end{equation}
Let us find the kernel of $p^1$. Let $\varphi : (\mathbb{R}^m,0) \to (\mathbb{R}^m,0)$ be an element of $\mathcal{D}_G$. 
Then, the map $g(u^k) = \|\frac{\partial \varphi^k}{\partial u^i}\|$ takes values in $G$, and $g(0) = I \in G$. Therefore, 
\begin{equation*}
\left.\frac{\partial^2\varphi}{\partial u^i\partial u^j}\right|_0 = \left.\frac{\partial g^k_i}{\partial u^j}\right|_0 
\end{equation*}
is a linear map $t : \mathbb{R}^m \cong T_0 \mathbb{R}^m \to \mathfrak{g}(G) \subset \mathfrak{gl}(m)$ 
but with property that $t(u)v = t(v)u$. In other words, elements of $\ker p^1$ are tensors of type $(2,1)$ on $\mathbb{R}^m$ such that $t^k_{ij} \in \mathfrak{g}$ for each $i$, and $t^k_{ij} = t^k_{ji}$.
The vector space of such tensors is called the \emph{first prolongation of the Lie algebra $\mathfrak{g}$} and is denoted by $\mathfrak{g}^1$.
From this follows that $\ker p^1$ is a commutative Lie group. 

Thus we have the exact sequence of Lie groups 
\begin{equation*}
0 \to \mathfrak{g}^1  \to G^1 \to G \to e, 
\label{eq:exact_splitting_sequence}
\end{equation*}
This sequence admits a splitting: for any $\|g^k_i\| \in G$ we take the diffeomorphism $\varphi^k(u^i) = g^k_i u^i$, which is evidently lies in $\mathcal{D}_G$,  and set $s(j^1_0\varphi) = j^2_0\varphi$. 

According to the group theory, the exact splitting sequence  \eqref{eq:exact_splitting_sequence} determines the right action of $G$ on $\mathfrak{g}^1$: $R_g t = s(g^{-1}) \cdot t \cdot s(g)$, and $G^1$ is the extension of $G$ by the commutative group $\mathfrak{g}^1$ with respect to the action $R$. 
This means that we have the group isomorphism
\begin{equation}
G^1 \to G \times \mathfrak{g}^1, \quad g^1 \mapsto (p^1(g^1), s\left( p^1( (g^1)^{-1})\right)g^1,
\label{eq:G_1_as_extension}
\end{equation}
therefore 
\begin{equation}
G^1 \cong G \times \mathfrak{g}^1, \text{ and } (g_1,t_1) \cdot (g_2,t_2) = (g_1g_2, R(g_2)t_1 + t_2).  
\label{eq:G_1_as_extension_1}
\end{equation}

Let us express this representation of $G^1$ in terms of the canonical coordinates. 
We have $p^1(\varphi^k_i,\varphi^k_{ij})=\varphi^k_i$, and $s(\varphi^k_i) = (\varphi^k_i,0)$. 
Hence follows that the isomorphism \eqref{eq:G_1_as_extension} is  
\begin{equation*}
(\varphi^k_i, \varphi^k_{ij}) \to (g^k_i, t^k_{ij})  \text{ with } 
g^k_i = \varphi^k_i, \ t^k_{ij} = \widetilde{\varphi}^k_s \varphi^s_{ij}.
\end{equation*}
Thus we get \emph{algebraic coordinates} $(g^k_i,t^k_{ij})$ on $G^1$ which are adopted to the algebraic structure of $G^1$. 

The right action $R$ with respect to the canonical coordinates is written as follows (we use \eqref{eq:product_D_2}):
\begin{equation*}
(\widetilde{\varphi}^k_i, 0) \cdot (\delta^k_i, \varphi^k_{ij}) \cdot (\varphi^k_i,0) = 
(\delta^k_i, \widetilde{\varphi}^k_s \psi^s_{lm} \varphi^l_i \varphi^m_j).
\end{equation*}
At the same time, the algebraic coordinates of the elements $(\varphi^k_i,0)$ and $(\delta^k_i,0)$ are the same, this means that they are $(g^k_i = \varphi^k_i, 0)$ and $(\delta^k_i, t^k_{ij} = \varphi^k_{ij})$.
Therefore, with respect to the algebraic coordinates the product of the group $G^1$ looks like
(see \eqref{eq:G_1_as_extension_1}):
\begin{equation}
(g^k_i, t^k_{ij}) \cdot (h^k_i, q^k_{ij}) = (g^k_s h^s_i, \widetilde{h}^k_s t^s_{pq} h^p_i h^q_j + q^k_{ij}).   
\label{eq:product_G_1_wrt_algebraic_coordinates}
\end{equation}

\subsection{Algebraic coordinates on $B^2(M)$. Description of $P^1(M)$ with respect to algebraic coordinates}
One can easily see that $D^2(m) = (GL(m))^1$. Therefore, we can consider the algebraic coordinates on $D^2(m)$ which 
give rise to the \emph{algebraic coordinates on the total space $B^2(M)$}.
Namely, if $(U,u)$ is a coordinate chart on $M$, and $(u^k, u^k_i, u^k_{ij})$ are the corresponding natural coordinates on $(\pi^2)^{-1}(U)$, for the \emph{algebraic coordinates} on $(\pi^2)^{-1}(U)$ we take
\begin{equation}
(u^k, p^k_i, p^k_{ij}) \text{ where } p^k_i = u^k_i, p^k_{ij} = \widetilde{u}^k_s u^s_{ij}. 
\label{eq:algebraic_coordinates_on_B_2}
\end{equation}
In fact, we change coordinates on the second factor of $U \times D^2(m) \cong (\pi^2)^{-1}(U)$.

With respect to the algebraic coordinates the first prolongation $P^1$ of an integrable $G$-structure $P$ is described in the following way: 
\begin{equation}
\left. P^{1} \right|_U = \left\{ (u^k, p^k_i, p^k_{ij}) \mid \|p^k_i\| \in G, \|p^k_{ij}\| \in \mathfrak{g}^1 \right\}.
\label{eq:P_1_wrt_algebraic_coordinates}
\end{equation}

If we have two coordinate charts of the atlas $\mathcal{A}$ on $M$ and $v^k = v^k(u^i)$ is the coordinate change, the corresponding gluing map \eqref{eq:gluing_maps_for_B_2} of trivializing charts of $P^1$  is written with respect to the algebraic coordinates as follows:
\begin{equation}
g : U \cap U' \to G^1, \quad 
g(x) = 
\left(\frac{\partial v^k}{\partial u^i}(u(x)), \frac{\partial u^k}{\partial v^s}(v(x)) \frac{\partial v^s}{\partial u^i \partial u^j}(u(x))\right).
\label{eq:gluing_map_of_P_2_wrt_algebraic_coordinates}
\end{equation}

\subsection{$P^1(M)$ in terms of $P(M)$}

\section{First prolongation of $G$-structure and associated bundles.}

Let us now consider the constructions of the previous section for cases $k=1$ and $k=2$. 

The case $k=1$ is rather simple.
We have $D^1(m) \cong GL(m)$ and $B^1(M)$ is the coframe bundle of $M$, 
because to each $j^1_x f$ we can put in correspondence the coframe $\{f^* du^i\}$ at $x \in M$. 

So we will consider in details the case $k=2$.  

\subsection{First prolongation of the Lie subgroup $G \subset GL(m)$}
\label{subsec:first_prolongation_of_lie_subgroup_GL(m)}

Now let us consider the set $\widetilde{D}^2(m) = (\varphi^k_i, \varphi^k_{ij})$, where $\|\varphi^k_i\|$ is an invertible matrix, and $\varphi^k_{ij}$ are not necessarily symmetric with respect to lower indices. 
Then, the set $\widetilde{D}$ endowed with the operation $*$ given by \eqref{eq:product_D_2} is a Lie group called the \emph{nonholonomic differential group of second order}.
We will call it \emph{first prolongation of the group} $GL(m)$ and denote by $GL^{(1)}(m)$.

On the group $GL^{(1)}(m)$ we can introduce another coordinate system:
\begin{equation}
g^k_j = \varphi^k_j, \quad a^k_{ij} = \widetilde{\varphi}^k_s \varphi^s_{ij} 
\label{eq:new_coorinates_differential_group}
\end{equation}
Then, with respect to these coordinates, by \eqref{eq:product_D_2}, we see  that the product $*$ can be written as follows:
\begin{equation}
(g^k_i, a^k_{ij}) * (h^k_i, b^k_{ij}) = (g^k_s h^s_i, \widetilde{h}^k_s a^s_{pq} h^p_i h^q_j + b^k_{ij}).   
\label{eq:product_new_coordinates}
\end{equation}
These formulas can be written in a matrix form. To do this, we consider $a^k_{ij}$ as a map $a : \mathbb{R}^m \to \mathfrak{gl}(m)$, $w^k \mapsto a^k_{ij}w^j$. 
Then, \eqref{eq:product_new_coordinates} takes the form
\begin{equation}
(g,a)*(h,b) = (g h, ad h^{-1} a \circ h + b).
\label{eq:product_matrix_form}
\end{equation}
Therefore, the group 
\begin{equation*}
GL^{(1)} \cong  (GL(m)\times Hom_{\mathbb{R}}(\mathbb{R}^m,\mathfrak{gl}(m)),*),
\end{equation*}
where $*$ is defined in \eqref{eq:product_matrix_form}.
\begin{remark}
The vector space $Hom_{\mathbb{R}}(\mathbb{R}^m,\mathfrak{gl}(m))$ is a right $GL(m)$-module with respect to the action 
\begin{equation}
\forall g \in GL(m), \quad a \in Hom_{\mathbb{R}}(\mathbb{R}^m,\mathfrak{gl}(m)), \quad g \cdot a = ad g^{-1} a g.
\label{eq:hom(R_m,gl(m))_is_a_right_GL(m)_module}
\end{equation}
The group $(GL(m)\times Hom_{\mathbb{R}}(\mathbb{R}^m,\mathfrak{gl}(m)),*)$ is the extension of the group $GL(m)$ 
with the right $GL(m)$-module $Hom_{\mathbb{R}}(\mathbb{R}^m,\mathfrak{gl}(m)),*)$. This fact motivates the definition of first prolongation for a subgroup $G \subset GL(m)$.
\end{remark}

\begin{remark}
The same considerations can be done for the holonomic jet group $D^2(m)$.
\end{remark}

\subsubsection{First prolongation of a Lie subgroup $G \subset GL(m)$}
The considerations of the previous subsection motivate 
\begin{definition}
Let $G \subset GL(m)$ be a Lie subgroup. Then the \emph{first prolongation of G} is the group 
\begin{equation}
G^{(1)} = G \times Hom_{\mathbb{R}}(\mathbb{R}^m,\mathfrak{g}(G))
\end{equation}
with product
\begin{equation*}
(g_1,a_1)(g_2,a_2) = (g_1 g_2, ad g_2^{-1} a_1 g_2 + a_2).
\end{equation*} 
\end{definition}
\begin{remark}
In this case also the vector space $Hom_{\mathbb{R}}(\mathbb{R}^m,\mathfrak{g}(G))$ is a right $G$-module with respect to the action 
\begin{equation}
\forall g \in G, \quad a \in Hom_{\mathbb{R}}(\mathbb{R}^m,\mathfrak{g}(G)), \quad g \cdot a = ad g^{-1} a g.
\label{eq:hom(R_m,g)_is_a_right_G_module}
\end{equation}
The group $(G\times Hom_{\mathbb{R}}(\mathbb{R}^m,\mathfrak{g}(G)),*)$ is the extension of the group $G$ 
with the right $G$-module $Hom_{\mathbb{R}}(\mathbb{R}^m,\mathfrak{g}(G)),*)$. 

Therefore we have the short exact sequence of Lie groups: 
\begin{equation*}
0 \to  Hom_{\mathbb{R}}(\mathbb{R}^m,\mathfrak{g}(G)) \to G^{(1)} \overset{\pi}{\longrightarrow} G.
\end{equation*}
The Lie homomorphism $\pi: G^{(1)} \rightarrow G$ comes from the natural projection of 2-jets onto 1-jets. 
\end{remark}

\subsubsection{Another coordinate system on $B^2(M)$}
Now recall that on $D^2(m)$ we can take another coordinate system $(g^k_i, a^k_{ij})$ 
(see \ref{eq:new_coorinates_differential_group}). 
With respect to the coordinates $(g^k_i, a^k_{ij})$, the $D^2(m)$-action on $B^2(M)$ is written as 
\begin{equation}
h^k_i = \widetilde{g}^k_s f^s_i,  \quad 
h^k_{ij} = \widetilde{g}^k_s f^s_{ij} - 
a^k_{lm} \widetilde{g}^l_t\widetilde{g}^m_r f^t_i f^r_j. 
\label{eq:action_of_D_2_wrt_coordinates_(g,a)}
\end{equation} 

Let us find the expression for a trivialization 
\begin{equation*}
(\pi^2)^{-1}(U) \to U \times D^2(m) 
\end{equation*}
with respect to the coordinates $(g,a)$.
For this purpose we take a section $s(x^i) = (x^i,\delta^k_i,0)$ of $B^2(M)$ over $U$, in fact this is the natural frame of order $2$, this means that it consists of $2$-jets of the coordinate functions. 
Then any point $b^2 = (x^i,f^k_i,f^k_{ij})$ can be written as $b^2 = s(x) \cdot (g,a)$, and the trivialization is given by
\begin{equation}
b^2 = (x^i,f^k_i,f^k_{ij}) \leftrightarrow (x, (g,a)).
\end{equation} 
Using \eqref{eq:action_of_D_2_wrt_coordinates_(g,a)}, and the coordinate expression $s(x)$, we get 
\begin{equation}
f^k_i = \widetilde{g}^k_i, \quad f^k_{ij} = - a^k_{lm} \widetilde{g}^l_i  \widetilde{g}^m_j. 
\end{equation}
Now let us write the gluing functions with respect to the coordinates $(p^k_i,p^k_{ij})$ on $B^2(M)$ 
and $(g^k_i,a^k_{ij})$. 
To do this we use the coordinate change  \eqref{eq:change_of_coordinates_f_to_p} on $B^2(M)$ and 
\eqref{eq:new_coorinates_differential_group} on the group $D^2(M)$.
We have
\begin{equation*}
\begin{split}
p^k_i = \widetilde{f}^k_i, \quad p^k_{ij} = - f^k_{lm} \widetilde{f}^l_i \widetilde{f}^m_j
\\
\bar x^k_i = g^k_i, \quad \bar x^k_{ij} = g^k_s a^s_{ij}  
\end{split}
\end{equation*}
and so,
\begin{equation}
f^k_i = \widetilde{p}^k_i, \quad f^k_{ij} = - p^k_{lm} \widetilde{p}^l_i \widetilde{p}^m_j
\label{eq:inverse_coordinate_transformation_B^2(M)}
\end{equation}

\begin{equation}
\begin{split}
& (x^k,f^k_i,f^k_{ij}) \to (x^k, (p^k_i,p^k_{ij}) ), \text{ where }
\\
&p^k_i = \widetilde{f}^k_i, \quad p^k_{ij} = - f^k_{lm} \widetilde{f}^l_i \widetilde{f}^m_j.
\end{split}
\label{eq:change_of_coordinates_f_to_p}
\end{equation} 
Using these formulas, and \eqref{eq:action_of_D_2_wrt_coordinates_(g,a)}, we can write the $D^2(m)$-action with respect to the coordinates $(x^i,p^k_i,p^k_{ij})$. 
If $(x^i,(p^k_i,p^k_{ij}))\cdot(g,a) = (x^i,(q^k_i,q^k_{ij}))$, then 
\begin{equation}
q^k_i = p^k_s g^s_i
\quad
q^k_{ij} = \widetilde{g}^k_s p^s_{lm} g^l_i g^m_j + a^k_{ij}.
\label{eq:action_of_D_2_wrt_coordinates_p_and_(g,a)}
\end{equation}
So, as it should be, at the second argument we get the product of elements of the group $D^2(m)$ (cf.  
\eqref{eq:product_new_coordinates}).
Hence, first of all from  
we get that
\begin{equation}
 \widetilde{\bar p}{}^k_s g^s_i = \widetilde{p}{}^k_i, \text{ and hence }  \bar p^k_i = \widetilde{g}^k_s p^s_i.
\label{eq:first_order_coordinate_change}
\end{equation}
Now, from 
with 
\eqref{eq:inverse_coordinate_transformation_B^2(M)}, we get
\begin{equation*}
- \bar p^k_{rt}\, \widetilde{\bar p}{}^r_l\widetilde{\bar p}{}^t_m g^l_i g^m_j + 
\widetilde{\bar p}{}^k_s g^s_t a^t_{ij} = f^k_{ij},
\end{equation*}
hence follows
\begin{equation*}
-\bar p^k_{lm} \widetilde{p}^l_i \widetilde{p}^m_j + \widetilde{p}^k_s a^s_{ij} = f^k_{ij} = 
- p^k_{lm} \widetilde{p}^l_i \widetilde{p}^m_j.
\end{equation*}
Finally, we obtain
\begin{equation*}
\bar p^k_{lm} = p^k_{lm} + \widetilde{p}^k_t a^t_{lm} p^l_i p^m_j.
\end{equation*}
As the result, we get the following theorem.

\begin{theorem}

\label{}
\end{theorem}

Note that 
\begin{equation}
x^k_i = \widetilde{\bar x^k_i}, 
\label{eq:correspondence_between_x_and_bar_x}
\end{equation}
With this notation, we have
\begin{equation}
\bar f^k_i = 
\label{eq:change_coordinates_B_2(M)}
\end{equation}

\section{Prolongation of $G$-structure}
\subsection{First prolongation of integrable $G$-structure}
Let $P(M,G)$ be an integrable $G$-structure, this means that there exists an atlas $(U_\alpha,u_\alpha)$ such that $\left\{ \frac{\partial}{\partial u_\alpha} \right\}$ are the sections of $P$. 

A first prolongation of $P(M,G)$ is the subbundle in $B^2(M)$ with the total space 
\begin{equation}
P^1(M) = \left\{ j^2_x f  \mid \left( \frac{\partial (f \circ u^{-1})^k}{\partial u^i}|_{u(x)} \right) \in G  \right\}
\end{equation}
and the structure group $G^1 \subset D^2(m)$ (the holonomic prolongation of $G$): 
\begin{equation}
G^1 = \left\{ j^2_0 \varphi \in D^2(m) \mid \left(\frac{\partial\varphi^k}{\partial u^i}\right) \in G \right\}.
\end{equation}

The product in $G^1$ is induced by the chain rule: if 
$j^k_0 \varphi = (\varphi^k_i, \varphi^k_{i_1i_2}, \dots, \varphi^k_{i_1,\dots,i_k})$
$j^k_0 \psi = (\psi^k_i, \psi^k_{i_1i_2}, \dots, \psi^k_{i_1,\dots,i_k})$,
then 
\begin{equation}
\eta^k_i = \psi^k_s \varphi^s_i, 
\quad
\eta^k_{ij} = \psi^k_{pq} \varphi^p_i \varphi^q_j + \psi^k_s \varphi^s_{ij}, \dots
\end{equation}

\subsection{Prolongation of $G$-structure}
\subsubsection{Nonholonomic jet bundle}
A $k$-th order nonholonomic jet is a ``nonsymmetric'' Teylor series and is defined by coordinates:
\begin{equation}
f^j_0 = \left( f^j, f^j_i, f^j_{i_1 i_2}, \dots, f^j_{i_1,i_2,\dots,i_k} \right),
\end{equation}
where $f^j_{i_1 \dots i_j}$ are not supposed to be symmetric with respect to lower indices.  

\subsection{Nonholonomic differential group}
The nonholonomic $k$-th order differential group is 

$D^{(k)}(m) = \{\varphi^k_0 = \left( \varphi^j, \varphi^j_i, \varphi^j_{i_1 i_2}, \dots, \varphi^j_{i_1,i_2,\dots,i_k} \right) \mid \det \left(\varphi^j_i\right) \ne 0\}$,

and the product is also ``induced by the chain rule'': 
if $\eta^k_0 = \psi^k_0 \cdot \varphi^k_0$ and  
$\varphi^k_0 = (\varphi^j_i, \varphi^j_{i_1i_2}, \dots, \varphi^j_{i_1,\dots,i_k})$
$\psi^k_0 = (\psi^j_i, \psi^j_{i_1i_2}, \dots, \psi^j_{i_1,\dots,i_k})$,
then 
$\eta^j_i = \psi^j_s \varphi^s_i$, $\eta^j_{i_1i_2} = \psi^j_{p_1p_2} \varphi^{p_1}_{i_1} \varphi^{p_2}_{i_2} + \psi^j_s \varphi^s_{i_1i_2}$, \dots

It is clear that $D^k(m)$ is a Lie subgroup of $D^{(k)}(m)$, so we have the left $D^k(m)$-action on $D^{(k)}(m)$. 

A \emph{nonholonomic $k$-th order jet bundle} is the locally trivial bundle over $M$ with the fiber $D^{(k)}(m)$ 
associated with the $D^k(m)$-principal bundle $B^k(M)$ with respect to the left $D^k(m)$-action on $D^{(k)}(m)$.

\subsection{First prolongation of arbitrary $G$-structure}
Let $P(M,G)$ be a $G$-structure.
A \emph{first prolongation of $P(M,G)$} is the subbundle in $B^{(2)}(M)$ with the total space 
\begin{equation}
P^{(1)}(M) = \left\{ f^2_0 \in B^{(2)}(M) \mid \left( f^j_i \right) \in G  \right\}
\end{equation}
and the structure group $G^{(1)} \subset D^{(2)}(m)$ (the nonholonomic prolongation of $G$): 
\begin{equation}
G^{(1)} = \left\{ \varphi^2_0 \in D^2(m) \mid \left(\varphi^k_i\right) \in G \right\}.
\end{equation}

\section{First prolongation of $G$-structure in terms of coframe bundle}

The first prolongation $P^{(1)}$ of $P$ can be expressed in terms of $P$ in the following ways: 
\begin{equation}
P^{(1)} = \{p^1 : \mathbb{R}^m \to T_p P | \theta_p \alpha = 1_{\mathbb{R}^m} \},
\label{eq:P_1_1}
\end{equation}
or
\begin{equation}
P^{(1)} = \{\omega : T_p P \to \mathfrak{g}  | \omega \sigma_p = 1_{\mathfrak{g}} \}
\label{eq:P_1_2}
\end{equation}
or
\begin{equation}
P^{(1)} = \left\{ H_p \mid H_p \oplus V_p = T_p P\right\}.
\label{eq:P_1_3}
\end{equation}
We will mainly use the first representation \eqref{eq:P_1_1}, but note that the third representation \eqref{eq:P_1_2} 
says that,  geometrically, $P^1$ consists of tangent subspaces transversal to vertical subspaces, i.\,e, of connections. 

The projection  $\pi^1_0 : P^1 \to P$,  is defined in terms of \eqref{eq:P_1_1} as follows $(p^1 : \mathbb{R}^m \to T_p P) \mapsto b$.

\begin{theorem}[Algebraic structure of $G^{(1)}$]
$G^{(1)}$ is isomorphic to the extension of $G$ via the $G$-module $\mathcal{L}(\mathbb{R}^n,\mathfrak{g})$:
\begin{equation}
G^{(1)} = G \times \mathcal{L}(\mathbb{R}^n,\mathfrak{g}), \quad  (g_1,a_1)*(g_2,a_2) = (g_1 g_2, ad g_2^{-1} a_1 + a_2).
\end{equation}
\end{theorem}

Action of $G^{(1)}$ on $P^{(1)}$ is described in the following way:
for $b^1 : \mathbb{R}^m \to T_b P$, and $g^1=(g,a) \in G^1$:
\begin{equation}
R^{(1)}_{g^1} b^1 = dR_g (b^1 \circ g) + \sigma_{pg} \circ a \circ g
\end{equation}

\section{First prolongation of equivariant map}

\subsection{First prolongation of a $G$-space}
Let $V$ be a manifold, then the first prolongation of $V$ is  
\begin{equation}
V^{(1)} = \left\{ v^1 :  \mathbb{R}^m \to T_v V \mid v \in V \right\},
\end{equation}
Let $\rho : G \times V \to V$ be a left action, then the first prolongation of $\rho$ is  
\begin{equation}
\rho^{(1)} : G^{(1)} \times V^{(1)} \to V^{(1)}, \quad  
\rho^{(1)}(g^1,v^1) = [dL_g \circ v^1 + \sigma_{gv}\circ A] \circ g^{-1}. 
\end{equation}

\subsection{Prolongation of equivariant map} 
$f : P \to V$ is an equivariant map.

Prolongation of $f$ is $f^{(1)} : P^{(1)} \to V^{(1)}$, $f^{(1)}(b^1)  = df_{\pi^1(b^1)} \circ b^1$.  

An equivariant map $f : P \to V$ determines a section $s : M \to E$, where $\pi_E : E \to M$ is a bundle with standard fiber $V$ associated with $P$.

$f^1 : P^1 \to V^1$ maps $b^1 \in P^1$ to the coordinates of $(s(\pi(b)),(\nabla s)(\pi(b)))$ with respect to $b$, where $\nabla$ corresponds to $b^1$.

\begin{example}[Simple example: vector field on $\mathbb{R}^n$]
$V$ is a vector field on $\mathbb{R}^m$.
The corresponding equivariant map is 
\begin{equation}
f : B(\mathbb{R}^m) \to \mathbb{R}^m, b = \{\eta^a\},  f(b) = \{\eta^a(V(\pi(b)))\}.
\end{equation}
Then the first prolongation of $f$ is defined as follows.  
If $b^{(1)} : e_i \in \mathbb{R}^m \to 
\frac{\partial}{\partial x^i} +  \Gamma^k_{ij} \frac{\partial}{\partial x^k_j}$,
then
\begin{equation}
f^{(1)} (b^{(1)}) = (V, \nabla(\omega) V) = (V^i, \partial_j V^i + \Gamma^k_{js} V^s).
\end{equation}

Action of $G^1$ is written as follows  
\begin{equation}
(g,a)(V,\nabla(\omega) V) = (\tilde g^k_i V^i, \tilde g^k_s \nabla_m V^s g^m_i  + \tilde g^k_s a^s_{mj} V^j g^m_i),  
\end{equation}

At $x_0$ such that  $V(x_0) = 0$, we get the action 
\begin{equation}
(g,a) (0,\nabla(\omega) V) = (0, \tilde g^k_s \nabla_m V^s g^m_i) = (0, \tilde g^k_s \partial_m V^s g^m_i)
\end{equation}

Therefore, the prolonged action coincides with the action of the group $GL(n)$ on the vector space of $n\times n$-matrices by conjugation. The invariants of this action are well known, for example, these are 
the trace and the determinant. 
Therefore, to find the invariants of a zero $x_0$ of a vector field $V$ we have to 
find the matrix $[\partial_i V^j(x_0)]$ and then write the invariants of this matrix under the conjugation, for example, one of them is $\det [\partial_i V^j(x_0)]$.
\end{example}

\end{document}